\nonstopmode
\documentclass[]{amsart}
\usepackage{amssymb}
\usepackage{amscd}
\usepackage{hyperref}
\usepackage{graphicx}
\usepackage{psfrag}
\usepackage{verbatim}
\usepackage[all]{xy}
\newdir{ >}{{}*!/-10pt/@{>}}
\newdir^{ (}{{}*!/-5pt/@^{(}}
\newdir_{ (}{{}*!/-5pt/@_{(}}

\theoremstyle{plain}
\newtheorem*{lemma*}{Lemma}
\newtheorem{lemma}[subsection]{Lemma}
\newtheorem*{theorem*}{Theorem}
\newtheorem{theorem}[subsection]{Theorem}
\newtheorem*{proposition*}{Proposition}

\newtheorem*{corollary*}{Corollary}
\newtheorem{corollary}[subsection]{Corollary}
\theoremstyle{definition}
\newtheorem*{definition*}{Definition}

\newtheorem*{example*}{Example}

\theoremstyle{remark}
\newtheorem*{remark*}{Remark}

\newenvironment{demo}[1]{\par\smallskip{\bf #1.}}{\par\smallskip}

\def\al{\alpha}
\def\be{\beta}

\def\et{\eta}
\def\vartheta{\theta}

\def\ka{\kappa}

\def\rh{\rho}

\def\ta{\tau}
\def\ph{\phi}

\def\om{\omega}

\def\Th{\Theta}

\def\i{^{-1}}
\def\x{\times}
\def\on{\operatorname}
\def\p{\partial}
\def\AMSonly#1{}
\begin{document}
\title{Homology and modular classes of Lie algebroids
}
\author{Janusz Grabowski, Giuseppe Marmo, Peter W. Michor  }
\address{
J. Grabowski: Mathematical Institute, Polish Academy of Sciences,
\'Sniadeckich 8, P.O. Box 21, 00-956 Warszawa, Poland }
\email{jagrab@impan.gov.pl}
\address{
G.\  Marmo: Dipartimento di Scienze Fisice, Universit\'a di Napoli
Federico II and INFN, Sezione di Napoli, via Cintia, 80126 Napoli, Italy }
\email{marmo@na.infn.it}
\address{
P.\  Michor: Institut f\"ur Mathematik, Universit\"at Wien, 
Nordbergstrasse~15, 
A-1090 Wien, Austria; {\it and:} Erwin Schr\"odinger Institut f\"ur
Mathematische Physik, Boltzmanngasse 9, A-1090 Wien, Austria }
\email{Peter.Michor@esi.ac.at }
\date{{\today} }
\thanks{The research of J.~Grabowski supported by the Polish Ministry of Scientific Research
and Information Technology under the grant No. 2 P03A 020 24. The research
of P.~Michor was supported by FWF Projects P~14195 and P 17108}

\begin{abstract} For a Lie algebroid, divergences chosen in a classical way lead
to a uniquely defined homology theory. They define also, in a natural way,
modular classes of certain Lie algebroid morphisms. This approach, applied
for the anchor map, recovers the concept of modular class due to S.~Evens,
J.-H.~Lu, and A.~Weinstein.
\end{abstract}
\def\LaTeXonly{}

\maketitle

\section{Introduction} 
Homology of a Lie algebroid structure on a vector bundle $E$ over $M$
are usually considered as homology of the corresponding Batalin-Vilkovisky
algebra associated with a chosen generating operator $\partial$ for the
Schouten-Nijenhuis bracket on multisections of $E$. The generating
operators that are homology operators, i.e. $\p^2=0$, can be identified
with flat $E$-connections on $\bigwedge^{\text{top}}E$ (see \cite{Xu}) or
divergence operators (flat right $E$-connections on $M\times\mathbb R$,
see \cite{Hu}). The problem is that the homology group depends on the
choice of the generating operator (flat connection, divergence) and no one
seems to be privileged. For instance, if a Lie algebroid on $T^*M$
associated with a Poisson tensor $P$ on $M$ is concerned, then the
traditional Poisson homology is defined in terms of the Koszul-Brylinski
homology operator $\p_P=[d,i_P]$. However, the Poisson homology groups may
differ from the homology groups obtained by means of 1-densities on $M$.
The celebrated modular class of the Poisson structure \cite{We} measures
this difference. Analogous statement is valid for triangular Lie
bialgebroids \cite{KS}.

The concept of a Lie algebroid divergence, so a generating operator,
associated with a `volume form', i.e. nowhere-vanishing section of
$\bigwedge^{\text{top}}E^*$, is completely classical (see \cite{KS},
\cite{Xu}). Less-known seems to be the fact that we can use `odd-forms'
instead of forms (cf. \cite{dR}) with same formulas for divergence and
that such nowhere-vanishing volume odd-forms always exist. The point is
that the homology groups obtained in this way are all isomorphic,
independently on the choice of the volume odd-form. This makes the
homology of a Lie algebroid a well-defined notion. From this point of view
the Poisson homology is not the homology of the associated Lie algebroid
$T^*M$ but a deformed version of the latter, exactly as the exterior
differential $d^\phi\mu=d\mu+\phi\wedge\mu$ of Witten \cite{Wi} is a
deformation of the standard de Rham differential.

In this language, the modular class of a Lie algebroid morphism
$\kappa:E_1\rightarrow E_2$ covering the identity on $M$ is defined as the
class of the difference between the pull-back of a divergence on $E_2$ and
a divergence on $E_1$, both associated with volume odd-forms. In the case
when $\kappa:E\rightarrow TM$ is the anchor map, we recognize the standard
modular class of a Lie algebroid \cite{ELW} but it is clear that other
(canonical) morphisms will lead to other (canonical) modular classes.

\section{Divergences and generating operators }\label{nmb:2}

\subsection{Lie algebroids and their cohomology}\label{nmb:2.1} 
Let $\ta:E\to M$ be a vector bundle. Let $\mathcal
A^i(E)=\on{Sec}(\bigwedge^iE)$ for $i=0,1,2,\dots$, let $A^i(E)=\{0\}$ for
$i<0$, and denote by $\mathcal A(E)=\bigoplus_{i\in \mathbb Z}\mathcal
A^i(E)$ the Grassmann algebra of multisections of $E$. It is a graded
commutative associative algebra with respect to the wedge product.

There are different ways to define a Lie algebroid structure on $E$. We
prefer to see it as a linear graded Poisson structure on $\mathcal A(E)$
(see \cite{GM1}), i.e., a graded bilinear operation $[\quad,\quad]$ on
$\mathcal A(E)$ of degree $-1$ with the following properties:
\begin{enumerate}
\item[(a)]
   Graded anticommutativity: $[a,b]=-(-1)^{(|a|-1)(|b|-1)}[b,a]$.
\item[(b)]
   The graded Jacobi identity:
   $[a,[b,c]]=[[a,b],c]+(-1)^{(|a|-1)(|b|-1)}[b,[a,c]]$.
\item[(c)]
   The graded Leibniz rule:
   $[a,b\wedge c]=[a,b]\wedge c + (-1)^{(|a|-1)|b|}b\wedge [a,c]$.
\end{enumerate}
This bracket is just the Schouten bracket associated with the the standard
Lie algebroid bracket on sections of $E$. It is well known that such
brackets are in bijective correspondence with de Rham differentials $d$ on
the Grassmann algebra $\mathcal A(E^*)$ of multisections of the dual
bundle $E^*$ which are described by the formula
\begin{multline} \label{fml:1}
d\mu(X_0,\dots,X_n) = \sum_i(-1)^i[X_i,\mu(X_0,\dots\hat{{}_i}\dots,
X_n)]+
\\
+\sum_{k<l}
(-1)^{k+l}\mu([X_k,X_l],X_0,\dots\hat{{}_k}\dots\hat{{}_l}\dots,X_n)
\end{multline}
where the $X_i$ are sections of $E$. We will refer to elements of
$\mathcal A(E^*)$ as {\it forms}. Since $d$ is a derivation on
$\mathcal A(E^*)$ of degree 1 with $d^2=0$, it defines the corresponding
de~Rham cohomology $H^*(E,d)$ of the Lie algebroid in the obvious way.

\subsection{Generating operators and divergences}\label{nmb:2.2} 
The definition of the homology of a Lie algebroid is more delicate than
that of cohomology. The standard approach is via generating operators for
the Schouten bracket $[\quad,\quad]$. By this we mean an operator $\p$ of
degree $-1$ on $\mathcal A(E)$ which satisfies
\begin{equation}\label{fml:2}
[a,b]=(-1)^{|a|}(\p(a\wedge b)-\p(a)\wedge b-(-1)^{|a|}a\wedge \p(b)).
\end{equation}
The idea of a generating operator goes back to the work by Koszul
\cite{Ko}. A generating operator which is a homology operator, i.e.
$\p^2=0$, gives rise to the so called Batalin-Vilkovisky algebra. Remark
that the leading sign $(-1)^{|a|}$ serves to produce graded antisymmetry
with respects to the degrees shifted by $-1$ out of graded symmetry. One
could equally well use $(-1)^{|b|}$ instead of $(-1)^{|a|}$, or one could
use the obstruction for $\p$ to be a graded right derivation in the
parentheses instead of a graded left one as we did. We shall stick to the
standard conventions.

It is clear from \thetag{\ref{fml:2}} and from the properties of the
Schouten bracket that $\p$ is then a second order differential operator on
the graded commutative associative algebra $\mathcal A(E)$ which is
completely determined by its restriction to $\on{Sec}(E)$. In fact, it is
easy to see (cf.\  \cite{Hu}) that
\begin{multline}
\label{fml:3} \p(X_1\wedge \dots\wedge X_n)
  = \sum_i(-1)^{i+1}\p(X_i)X_1\wedge\dots\hat{{}_i}\dots\wedge X_n +
\\
+\sum_{k<l}(-1)^{k+l}[X_k,X_l]\wedge
X_1\wedge\dots\hat{{}_k}\dots\hat{{}_l}\dots\wedge X_n
\end{multline}
for $X_1,\dots,X_n\in \on{Sec}(E)$, which looks completely dual to
\thetag{\ref{fml:1}}. From \thetag{\ref{fml:2}} we get the following
property of $\p$:
\begin{equation}
\label{fml:4} -\p(fX) = -f\p(X)+[X,f]\quad\text{  for }X\in\on{Sec}(E),
f\in C^\infty(M).
\end{equation}
Since $[X,f]=\rh(X)(f)$ where $\rh:E\to TM$ is the anchor map of the Lie
algebroid structure on $E$, the operator $-\p$ has the algebraic property
of a divergence. Conversely, \thetag{\ref{fml:3}} defines a generating
operator for $[\cdot,\cdot]$ if only \thetag{\ref{fml:4}} is satisfied.
i.e., generating operators can be identified with divergences. We may
express this by $\on{div}\leftrightarrow \p_{\on{div}}$. But a true
divergence $\on{div}:\on{Sec}(E)\to C^\infty(M)$ satisfies besides
\thetag{\ref{fml:4}} a cocycle condition
\begin{equation}\label{fml:5}
\on{div}([X,Y]) = [\on{div}(X),Y] + [X,\on{div}(Y)], \quad
X,Y\in\on{Sec}(E),
\end{equation}
which is equivalent (see \cite{Hu}) to the fact that the corresponding
generating operator $\p_{\on{div}}$ is a homology operator:
$(\p_{\on{div}})^2=0$. Note that divergences can be used in construction
of generating operators also in the supersymmetric case (cf. \cite{KSMo}).

From now on we will fix the Lie algebroid structure on $E$, and we will
denote by $\on{Gen}(E)$ the set of generating operators for
$[\cdot,\cdot]$ which are homology operators, and by $\on{Div}(E)$ the
canonically isomorphic (by \thetag{\ref{fml:3}}) set of divergences for
the Lie algebroid satisfying \thetag{\ref{fml:4}} and
\thetag{\ref{fml:5}}. The problem is that there does not exist a canonical
divergence, thus no canonical generating operator.

The set $\on{Div}(E)$ can be identified with the set of all flat
$E$-connections on $\bigwedge^{\text{top}}E^*$, i.e., operators
$\nabla:\on{Sec}(E)\x
\on{Sec}(\bigwedge^{\text{top}}(E^*))\rightarrow\on{Sec}
(\bigwedge^{\text{top}}(E^*))$ which satisfy \begin{enumerate}
\item[(i)]
$\nabla_{fX}\mu = f\nabla_X\mu$,
\item[(ii)]
$\nabla_{X}(f\mu) = f\nabla_X\mu + \rh(X)(f)\mu$,
\item[(iii)]
$[\nabla_X,\nabla_Y] = \nabla_{[X,Y]}$.
\end{enumerate}
The identification is via
\begin{equation}\label{fml:6}
\mathcal L_X\mu-\nabla_X\mu = \on{div}(X)\mu
\end{equation}
(cf. \cite[(50)]{KS}), where $\mathcal L_X=d\,i_X+i_X\,d$ is the Lie
derivative. Note that \thetag{\ref{fml:6}} is independent of the choice of
the section $\mu\in\on{Sec}(\bigwedge^{\text{top}}(E^*))$. We can use
$\bigwedge^{\text{top}}(E)$ instead of $\bigwedge^{\text{top}}(E^*)$ and
get the identification of $\on{Div}(E)$ with the set of flat
$E$-connection on $\bigwedge^{\text{top}}(E)$ by (see \cite{Xu})
\begin{equation}\label{fml:6a}
\mathcal L_X\Lambda-\nabla_X\Lambda = \on{div}(X)\Lambda.
\end{equation}

Of course, additional structures on $E$ as, e.g., a Riemannian metric
(smoothly arranged scalar products on fibers of $E$), may furnish a
distinguished divergence on $E$. Fixing a metric we can distinguish a
canonical torsionfree connection $\nabla$ on $E$ - the Levi-Civita
connection for the Lie algebroid - in the standard way. It satisfies the
standard Bianchi and Ricci identities (see \cite{Ne}) and induces a
connection on $\bigwedge^{\text{top}}(E)$ for which the generating
operator $\p_\nabla$ has the local form (see \cite{Xu})
$\p_\nabla(a)=-\sum_k i(\al^k)\nabla_{X_k}a$ where the $X_k$ and $\al^k$
are dual local frames for $E$ and $E^*$, respectively. Since
$$\p_\nabla^2=\sum_{k,j}i(\al^j)\nabla_{X_j}i(\al^k)\nabla_{X_k}=
\sum_{k,j}i(\al^j)i(\al^k)(\nabla_{X_j}\nabla_{X_k}-\nabla_{\nabla_{X_j}X_k}),
$$
$\p^2=0$ is equivalent to
\begin{equation}\label{fml:7}
\sum_{j,k}i(\al^j)i(\al^k)R(X_j,X_k)=0,
\end{equation}
where $R$ is the curvature tensor of $\nabla$. For a Levi-Civita
connection $\nabla$ the generating operator $\p_\nabla$ is really a
homology operator due to the following lemma.

\begin{lemma}\label{nmb:2.3}
A torsionfree connection $\nabla$ on $E$ satisfies simultaneously the
Bian\-chi and the Ricci identity if and only if  \thetag{\ref{fml:7}}
holds for dual local frames $X_k$ and $\al^k$ of $E$ and $E^*$,
respectively.
\end{lemma}

\begin{demo}{Proof}
\thetag{\ref{fml:7}} is equivalent to $\sum_{j,k}R(X_j,X_k)^*(\al^k\wedge
\al^j\wedge \om)=0$ for all forms $\om$. It suffices to check this for
$\om$ a function or a 1-form due to the derivation property of
contractions. For $\om$ a function $f$ we have
\begin{align*}
\sum_{j,k}R(X_j,X_k)^*(f\al^k\wedge \al^j) &=
\sum_{j,k}f\Bigl(R(X_j,X_k)^*(\al^k)\wedge \al^j
  +\al^k\wedge R(X_j,X_k)^*(\al^j)\Bigr)
\\&
= 2f\sum_{s,j,k}R^k_{jks}\al^s\wedge \al^j
\end{align*}
and this vanishes for all $f$ if and only if $R^k_{jks}$ is symmetric in
$(j,s)$, i.e., if the Ricci identity holds. For $\om$ a 1-form, say
$\al^i$, we have
\begin{multline*}
\sum_{j,k}R(X_j,X_k)^*(\al^k\wedge \al^j\wedge \al^i)
=\\
= \sum_{j,k}\Bigl(R(X_j,X_k)^*(\al^k\wedge \al^j)\wedge \al^i
  +\al^k\wedge \al^j\wedge R(X_j,X_k)^*(\al^i)\Bigr)
=\\
= 0+ \sum_{j,k,s}R^i_{jks}\al^k\wedge \al^j \wedge \al^s
\end{multline*}
and this vanishes for all $i$ if and only if
$\sum_{\text{cycl}(j,k,s)}R^i_{jks}=0$, i.e., if the first Bianchi
identity holds. \qed\end{demo}

\begin{corollary}\label{nmb:2.4}
Any Levi-Civita connection for a Riemannian metric on a Lie algebroid $E$
induces a flat connection on $\bigwedge^{\text{top}}E$, thus also on
$\bigwedge^{\text{top}}E^*$. \qed
\end{corollary}

\section{Homology of the Lie algebroid }\label{nmb:3}

\subsection{Getting divergences from odd forms}\label{nmb:3.1} 
There is no distinguished divergence for the Lie algebroid structure on
$E$, but there is a distinguished subset of divergences which we may
obtain in a classical way. Firstly, suppose that the line bundle
$\bigwedge^{\text{top}}E^*$ is trivializable. So we can choose a vector
volume, i.e., a nowhere vanishing section
$\mu\in\on{Sec}(\bigwedge^{\text{top}}E^*)$. Then the formula
\begin{equation}\label{fml:8}
\mathcal L_X\mu = \on{div}_\mu(X)\mu,\quad\text{ where }X\in\on{Sec}(E)
\end{equation}
defines a divergence $\on{div}_\mu$. We observe that
$\on{div}_{-\mu}=\on{div}_{\mu}$. Thus for the non-orientable case we look
for sections of a bundle over $M$ which locally consists of non-ordered
pairs $\{\mu_{\al},-\mu_{\al}\}$ for an open cover $M=\bigcup_\al U_{\al}$
such that the sets $\{\mu_{\al},-\mu_{\al}\}$ and
$\{\mu_{\be},-\mu_{\be}\}$ coincide when restricted to $U_\al\cap U_\be$.
The fundamental observation is that such global sections always exist and
define global divergences. This is because they can be viewed as sections
of the bundle $|\on{Vol}|_E=(\bigwedge^{\text{top}}E^*)_0/\mathbb Z_2$,
where $(\bigwedge^{\text{top}}E^*)_0$ is the bundle
$\bigwedge^{\text{top}}E^*$ with the zero section removed and divided by
the obvious $\mathbb Z_2$-action of passing to the opposite vector. The
bundle $|\on{Vol}|_E$ is an 1-dimensional affine bundle modelled on the
vector bundle $M\x \mathbb R$, and also a principal $\mathbb R$ bundle
where $t\in \mathbb R$ acts by scalar multiplication with $e^t$. Since it
has a contractible fiber, sections always exist. Note that sections
$|\mu|$ of $|\on{Vol}|_E$ are particular cases of odd forms, \cite{dR}:
Let $p:\tilde M\to M$ be the two-fold covering of $M$ on which $p^*E$ 
is oriented, namely the set of vectors of length 1 in the line bundle over 
$M$ with cocycle of transition functions $\on{sign}\det(\ph_{\al\be})$, where 
$\ph_{\al\be}:U_\al\cap U_\be\to GL(V)$ is the cocycle of transition 
functions for the vector bundle $E$. Then the odd forms are those forms 
on $p^*E$ which are in the $-1$ eigenspace of the natural vector bundle 
isomorphism which covers the decktransformation of $\tilde M$.
So odd forms are certain sections of a line bundle over a two-fold covering of 
the base manifold $M$. 
This is related but complementary to the construction of the line bundle 
(over $M$) of densities 
which involve the cocycle of transition functions $|\det(\ph_{\al\be})|$.
For example, any Riemannian metric $g$ on the vector bundle $E$
induces an odd volume form $|\mu|_g\in Sec(|\on{Vol}|_E)\simeq
Sec(|\on{Vol}|_{E^*})$ which locally is represented by the wedge product
of any orthonormal basis of local sections of $E$ (thus $E^*$). Note that
such product is independent on the choice of the basis modulo sign, so our
odd volume is well defined.

For the definition of a divergence $\on{div}_{|\mu|}$ associated to
$|\mu|\in\on{Sec}(|\on{Vol}|_E)$ we will write simply
\begin{equation}\label{fml:9}
\mathcal L_X|\mu|=\on{div}_{|\mu|}(X)|\mu| \quad\text{ for
}X\in\on{Sec}(E).
\end{equation}
Note that the distinguished set $\on{Div}_0(E)$ of divergences obtained in
this way from sections of $|\on{Vol}|_E$ corresponds (in the sense of
\thetag{\ref{fml:6}}) to the set of those flat connections on
$\bigwedge^{\text{top}}E^*$ whose holonomy group equals $\mathbb Z_2$:
Associate the horizontal leaf $|\mu|$ to such a connection, and note that
a positive multiple of $|\mu|$ gives rise to the same divergence.

In the case of a vector bundle Riemannian metric $g$ on $E$ a natural
question arises about the relation between the divergence
$\on{div}_{|\mu|_g}$ associated with the odd volume $|\mu|_g$ induced by
the metric $g$ and the divergence $\on{div}_{\nabla_g}$ induced by the
flat Levi-Civita connection $\nabla_g$ on $\bigwedge^{\text{top}}E^*\simeq
\bigwedge^{\text{top}}E$.

\begin{theorem}\label{nmb:3.2} 
For any vector bundle Riemannian
metric $g$ on $E$
$$\on{div}_{|\mu|_g}=\on{div}_{\nabla_g}.$$
\end{theorem}
\begin{demo}{Proof}
Let $X_1,\dots,X_n$ be an orthonormal basis of local sections of $E$ and
$\alpha^k=g(X_k,\cdot)$ be the dual basis of local sections of $E^*$, so
that $|\mu|_g$ is locally represented by $\alpha_1\wedge\dots\alpha^n$.
For any local section $X$ of $E$
\begin{multline*}\on{div}_{|\mu|_g}(X)=-\langle\mathcal
L_X(\alpha^1\wedge\cdots\wedge\alpha^n),X_1\wedge\cdots\wedge
X_n\rangle=\\
\langle\alpha^1\wedge\cdots\wedge\alpha^n,\mathcal
L_X(X_1\wedge\cdots\wedge X_n)\rangle=
\sum_k\langle\alpha^k,[X,X_k]\rangle=\\
\sum_k\langle\alpha^k,\nabla_XX_k-\nabla_{X_k}X\rangle=
\sum_kg(X_k,\nabla_XX_k)-\sum_ki(\alpha^k)\nabla_{X_k}X.
\end{multline*}
But $-\sum_ki(\alpha^k)\nabla_{X_k}X=\on{div}_{\nabla_g}(X)$ and
$$2\sum_kg(X_k,\nabla_XX_k)=\sum_k\rho(X)g(X_k,X_k)-\sum_k\nabla_X(g)(X_k,X_k)=0,$$
where $\rho:E\rightarrow TM$ is the anchor of the Lie algebroid on $E$,
since $\nabla$ is Levi-Civity ($\nabla g=0$). \qed\end{demo}

\subsection{The generating operator for an odd form}\label{nmb:3.4} 
The corresponding generating operator $\p_{|\mu|}$ for the divergence of
a non-vanishing odd form $|\mu|$ can be defined explicitly by
\begin{displaymath}
\mathcal L_a|\mu|=-i(\p_{|\mu|}(a))|\mu|,
\end{displaymath}
where $\mathcal L_a = i_a\, d - (-1)^{|a|}d\,i_a$ is the Lie differential
associated with $a\in\mathcal A^{|a|}(E)$ so that
\begin{equation}\label{fml:11}
i(\p_{|\mu|}(a))|\mu| = (-1)^{|a|}d\,i_a|\mu|.
\end{equation}
In other words, locally over $U$ we have
\begin{equation}\label{fml:12}
\p_{|\mu|}(a) = (-1)^{|a|} *_{\mu}\i\; d\; *_{\mu}(a)
\end{equation}
where $*_\mu$ is the isomorphism of $\mathcal A(E)|_U$ and $\mathcal
A(E^*)|_U$ given by $*_\mu(a)=i_a\mu$, for a representative $\mu$ of
$|\mu|$. Note that the right hand side of  \thetag{\ref{fml:12}} depends
only on $|\mu|$ and not on the choice of the representative since
$*_\mu\;d\;*_\mu=*_{-\mu}\;d\;*_{-\mu}$. Formula  \thetag{\ref{fml:12}}
gives immediately $\p_{|\mu|}^2=0$, which also follows from the remark on
flat connections above. So $\p_{|\mu|}$ is a homology operator.

Moreover, it is also a generating operator. Namely, using standard
calculus of Lie derivatives we get
\begin{displaymath}
\mathcal L_{a\wedge b} = i_b\;\mathcal L_a -(-1)^{|a|} i_{[a,b]} +
(-1)^{|a|\,|b|}i_a\;\mathcal L_b
\end{displaymath}
which can be rewritten in the form
\begin{equation}\label{fml:13}
i_{[a,b]}=(-1)^{|a|}\Bigl(-\mathcal L_{a\wedge b} +i_b\;\mathcal
L_a+(-1)^{|a|(|b|+1)}i_a\;\mathcal L_b\Bigr)
\end{equation}
When we apply  \thetag{\ref{fml:13}} to $|\mu|$ we get
\begin{displaymath}
i_{[a,b]}|\mu|=(-1)^{|a|}\Bigl(i(\p_{|\mu|}(a\wedge b))
-i(\p_{|\mu|}(a)\wedge b)-(-1)^{|a|}i(a\wedge \p_{|\mu|}(b))\Bigr)|\mu|
\end{displaymath}
which proves  \thetag{\ref{fml:2}}. Thus we get:

\begin{theorem}\label{nmb:3.3}
For any $|\mu|\in\on{Sec}(|\on{Vol}|_E)$ the formula
\begin{equation}\label{fml:10}
\mathcal L_a|\mu| = - i(\p_{|\mu|}(a))|\mu|
\end{equation}
defines uniquely a generating operator $\p_{|\mu|}\in \on{Gen}(E)$.
\end{theorem}

We remark that formula  \thetag{\ref{fml:10}} in the case of trivializable
$\bigwedge^{\text{top}}E^*$ has been already found in \cite{KS}. In this
sense the formula is well known. What is stated in Theorem \ref{nmb:3.3} is that
\thetag{\ref{fml:10}} serves in general, as if the bundle
$\bigwedge^{\text{top}}E^*$ were trivial, if we replace ordinary forms
with odd volume forms.

\subsection{Homology of the Lie algebroid}\label{nmb:3.5} 
The homology operator of the form $\p_{|\mu|}$ will be called the
homology operator for the Lie algebroid $E$. The crucial point is that
they all define the same homology. This is due to the fact that
$\p_{|\mu_1|}$ and $\p_{|\mu_2|}$ differ by contraction with an exact
1-form.

In general, two divergences differ by contraction with a closed 1-form.
Indeed, $(\on{div}_1-\on{div}_2)(fX)=f(\on{div}_1-\on{div}_2)(X)$, so
$(\on{div}_1-\on{div}_2)(X)=i_\ph X$ for a unique 1-form $\ph$. Moreover,
\thetag{\ref{fml:5}} implies that $i_\ph[X,Y]=[i_\ph X,Y]+[X,i_\ph Y]$, so
$\ph$ is closed. Since both sides are derivations we have
\begin{equation}\label{fml:14}
\p_{\on{div}_2}-\p_{\on{div}_1} = i_\ph.
\end{equation}
But for any $|\mu_1|,|\mu_2|\in\on{Sec}(|\on{Vol}|_E)$ there exists a
positive function $F=e^f$ such that $|\mu_2|=F\,|\mu_1|$. Then $\mathcal
L_X|\mu_2|=\mathcal L_X(F\,|\mu_1|)=\mathcal L_X(F)\,|\mu_1|+F\,\mathcal
L_X(|\mu_1|)$ so that $\on{div}_{|\mu_2|}(X)|\mu_2|=\mathcal
L_X(f)\,|\mu_2|+\on{div}_{|\mu_1|}(X)|\mu_2|$, i.e.,
$\on{div}_{|\mu_2|}-\on{div}_{|\mu_1|}=i(df)$.

To see that the homology of $\p_{|\mu_1|}$ and $\p_{|\mu_2|}$ are the
same, note first that $\p_{|\mu_2|} = \p_{|\mu_1|}a + i_{df}a$. And then
let us gauge $\mathcal A(E)$ by multiplication with $F=e^f$. This is an
isomorphism of graded vector spaces and we have
\begin{displaymath}
e^f\;\p_{|\mu_1|}\;e^{-f}\,a = \p_{|\mu_1|}a + i_{df}a=\p_{|\mu_2|}a,
\end{displaymath}
so $\p_{|\mu_1|}$ and $\p_{|\mu_2|}$ are graded conjugate operators.

This is just the dual picture of the well-known gauging of the de~Rham
differential by Witten \cite{Wi}, see also \cite{GM1} for consequences in
the theory of Lie algebroids. Thus we have proved (cf. \cite[p.120]{KS}):

\begin{theorem}\label{nmb:3.6} 
All homology operators for a Lie algebroid
generate the the same homology: $H_*(E,\p_{|\mu_1|}) =
H_*(E,\p_{|\mu_2|})$. In the case of trivializable
$\bigwedge^{\text{top}}E^*$,  \thetag{\ref{fml:12}} gives Poincar\'e
duality $H^*(E,d)\cong H_{\text{top}-*}(E,\p_{|\mu|})$.
\end{theorem}

\subsection{Remark}\label{nmb:3.7} 
We got a well-defined Lie algebroid homology, in contrast with the
standard approach when all generating operators are admitted. It is clear
that adding a term $i_\ph$ with $\ph$ a closed 1-from which is not exact,
as in \thetag{\ref{fml:14}}, will probably change the homology. But this
could be understood as an a priori deformation like in the case of the
deformed de~Rham differential of Witten \cite{Wi}:
\begin{equation}\label{fml:15}
d^\ph\et = d\et +\ph\wedge \et.
\end{equation}
Indeed, $i(i_\ph a)\mu=-(-1)^{|a|}\ph\wedge i_a\mu$ implies
$*_\mu\,i_\ph(a)=-(-1)^{|a|}\,e_\ph\,*_\mu(a)$ where $e_\ph\et=\ph\wedge
\et$. Thus we get
$(-1)^{|a|}\,*_\mu\i(d+e_\ph)\,*_\mu(a)=(\p_{|\mu|}-i_\ph)(a)$, so, at
least in the the trivializable case, there is the Poincar\'e duality
\begin{displaymath}
H^*(E,d+e_\ph)\cong H_{\text{top}-*}(E,\p_{\mu}-i_\ph).
\end{displaymath}
Note that the differentials $d^\ph$ appear as part of the Cartan
differential calculus for Jacobi algebroids, see \cite{IM}, \cite{GM},
\cite{GM1}, so that there is a relation between generating operators for a
Lie algebroid and the Jacobi algebroid structures associated with it.

\section{Modular classes}\label{nmb:4}

\subsection{The modular class of a morphism}\label{nmb:4.1} 
As we have shown, every Lie algebroid $E$ has a distinguished class
$\on{Div}_0(E)$ of divergences obtained from sections of $|\on{Vol}|_E$.
Such divergences differ by contraction with an exact 1-form. Let now
$\ka:E_1\to E_2$ be a morphism of Lie algebroids.

There is the induced map $\ka^*:\on{Div}(E_2)\to\on{Div}(E_1)$ defined by
$\ka^*(\on{div}_2)(X_1)=\on{div}_2(\ka(X_1))$. The fact that $\ka^*$ maps
divergences into divergences follows from $\ka(fX)=f\ka(X)$ and the fact
that the Lie algebroid morphism respects the anchors,
$\rho_1=\rho_2\circ\ka$. The space $\ka^*(\on{Div}_0(E_2))\subset
\on{Div}(E_1)$ consists of divergences which differ by insertion of an
exact 1-form. Therefore, the cohomology class of the 1-form $\ph$ which is
defined by the equation
\begin{equation}\label{fml:16}
\ka^*(\on{div}_{E_2})-\on{div}_{E_1} = i_\ph,\quad\text{  for
}\on{div}_{E_i}\in\on{Div}_0(E_i), i=1,2,
\end{equation}
does not depend on the choice of $\on{div}_{E_1}$ and $\on{div}_{E_2}$. We
will call it the {\it modular class} of $\ka$ and denote it by
$\on{Mod}(\ka)$. Thus we have:

\begin{theorem} 
For every Lie algebroid morphism
\begin{displaymath}
\xymatrix{
E_1 \ar[-0,2]^{\ka} \ar[1,1]^{\ta_1} & & E_2 \ar[1,-1]_{\ta_2} \\
& M & }
\end{displaymath} 
the cohomology class $\on{Mod}(\ka)=[\ph]\in
H^1(E_1,d_{E_1})$ defined by $\ph$ in  \thetag{\ref{fml:16}} is well
defined independently of the choice of $\on{div}_{E_1}\in\on{Div}_0(E_1)$
and $\on{div}_{E_2}\in\on{Div}_0(E_2)$. \qed
\end{theorem}

\subsection{The modular class of a Lie algebroid}\label{nmb:4.2} 
In the case when the morphism $\ka=\rh:E\to TM$ is the anchor map of a
Lie algebroid $E$, the modular class $\on{Mod}(\rh)$ is called the
{\it modular class of the Lie algebroid} $E$ and it is denoted by
$\on{Mod}(E)$. The idea that the modular class is associated with the
difference between the Lie derivative action on
$\bigwedge^{\text{top}}(E^*)$ and on $\bigwedge^{\text{top}}T^*M$  via the
anchor map is, in fact, already present in \cite{ELW}. Also the
interpretation of the modular class as certain secondary characteristic
class of a Lie algebroid, present in \cite{Fe}, is a quite similar. In
\cite{Fe} the trace of the difference of some connections is used instead
of the difference of two divergences. We have

\begin{theorem} $\on{Mod}(E)$ is the modular class $\Th_E$ in the sense
of \cite{ELW}.
\end{theorem}

\begin{demo}{Proof} The modular class $\Th_E$ in the sense of \cite{ELW} is
defined as the class $[\ph]$ where $\ph$ is given by
\begin{equation}\label{fml:17}
\mathcal L_X(a)\otimes\mu + a\otimes\mathcal L_{\rh(X)}\mu=\langle X,\ph
\rangle a\otimes \mu
\end{equation}
for all sections $a$ of $\bigwedge^{\text{top}}(E)$ and $\mu$ of
$\bigwedge^{\text{top}}(T^*M)$, respectively. Let us take
$|a^*|\in\on{Sec}(|\on{Vol}|_E)$ and $|\mu|\in\on{Sec}(|\on{Vol}|_{TM})$,
locally represented by $a^*\in\on{Sec}(\bigwedge^{\text{top}}(E^*|_U))$
and $\mu\in\on{Sec}(\bigwedge^{\text{top}}(T^*M|_U))$. Let $a$ be a local
section of $\bigwedge^{\text{top}}E$ dual to $a^*$. Then $\mathcal
L_X(a)=-\on{div}_{|a^*|}(X)\,a$ and $\mathcal
L_X(\mu)=\rh^*(\on{div}_{|\mu|})(X)\,\mu$ so that \thetag{\ref{fml:17}}
yields $i_\ph=\rh^*(\on{div}_{|\mu|})-\on{div}_{|a^*|}$. \qed\end{demo}

Note that in our approach the modular class $\on{Mod}(TM)$ of the
canonical Lie algebroid $TM$ is trivial by definition. It is easy to see
that the modular class of a base preserving morphism can be expressed in
terms of the modular classes of the corresponding Lie algebroids.

\begin{theorem} For a base preserving morphism
$\ka:E_1\rightarrow E_2$ of Lie algebroids
\begin{displaymath}
\on{Mod}(\ka)=\on{Mod}(E_1)-\ka^*(\on{Mod}(E_2)).\end{displaymath}
\end{theorem}

\begin{demo}{Proof} Let $\rho_l:E_l\rightarrow TM$ be the anchor of $E_l$,
$l=1,2$. Take $\on{div}_{E_l}\in\on{Div}_0(E_l)$, $l=1,2$, and
$\on{div}_{TM}\in\on{Div}_0(TM)$. Since $\on{Mod}(E_l)$ is represented by
$\eta_l$, $i_{\eta_l}=\on{div}_{E_l}-\rho_l^*(\on{div}_{TM})$ and
$\rho_1=\rho_2\circ\ka$, we can write
\begin{multline*}
i_{\eta_1}=\on{div}_{E_1}-\rho_1^*(\on{div}_{TM})=\\
\on{div}_{E_1}-\ka^*(\on{div}_{E_2})+\ka^*(\on{div}_{E_2}-\rho_2^*(\on{div}_{TM})=
i_{\eta_\ka}+i_{\ka^*(\eta_2)},
\end{multline*}
where $\eta_\ka$ represents $\on{Mod}(\ka)$. Thus
$\eta_1=\eta_\ka+\eta_2$. \qed\end{demo}

\subsection{The universal Lie algebroid}\label{nmb:4.3} 
For any vector bundle $\ta:E\to M$ there exists a universal Lie
algebroid $\on{QD}(E)$ whose sections are the quasi-derivations on $E$,
i.e., mappings $D:\on{Sec}(E)\to \on{Sec}(E)$ such that
$D(fX)=f\,D(X)+\hat D(f)\,X$ for $f\in C^\infty(M)$ and $X\in\on{Sec}(E)$,
where $\hat D$ is a vector field on $M$; see the survey article \cite{Gr}.
Quasi-derivations are known in the literature under various names:
covariant differential operators \cite{Ma}, module derivations \cite{Ne},
derivative endomorphisms \cite{KSMa}, etc. The Lie algebroid $\on{QD}(E)$
can be described as the Atiyah algebroid associated with the principal
$GL(n,\mathbb R)$-bundle $\on{Fr}(E)$ of frames in $E$, and
quasi-derivations can be identified with the $GL(n,\mathbb R)$-invariant
vector fields on $\on{Fr}(E)$. The corresponding short exact Atiyah
sequence in this case is
\begin{displaymath}
0\to \on{End}(E)\to \on{QD}(E) \to TM \to 0.
\end{displaymath}
This observation shows that there is a modular class associated to every
vector bundle $E$, namely the modular class $\on{Mod}(\on{QD}(E))$, which
is a vector bundle invariant.

It is also obvious that, viewing a flat $E_0$-connection (representation)
in a vector bundle $E$ over $M$ for a Lie algebroid $E_0$ over $M$ as a
Lie algebroid morphism $\nabla:E_0\rightarrow QD(E)$, one can define the
modular class $\on{Mod}(\nabla)$.

\medskip\noindent
\textbf{Question.} \textit{How is $\on{Mod}(\on{QD}(E))$ related to other
invariants of $E$ (e.g. characteristic classes)?}

\subsection{Remark}\label{nmb:4.4} 
One can interpret the modular class $\on{Mod}(E)$ of the Lie algebroid
$E$ as a "trace" of the adjoint representation. Indeed, if we fix local
coordinates $u^a$ on $U\subset M$ a local frame $X_i$ of local sections of
$E$ over $U$, and the dual frame $\al^i$ of $E^*$, then the Lie algebroid
structure is encoded in the "structure functions"
\begin{displaymath}
[X_i,X_j]= \sum_k c^k_{ij}\,X_k, \qquad \rh(X_i) = \sum_a
\rh^a_i\;\p_{u^a}.
\end{displaymath}

\begin{proposition*} The modular class $\on{Mod}(E)$ is locally
represented by the closed 1-form
\begin{equation}\label{fml:18}
\ph=\sum_{i}\Bigl(\sum_k c^k_{ik}+\sum_a\frac{\p \rh^a_i}{\p
u^a}\Bigr)\al^i.
\end{equation}
\end{proposition*}

\begin{demo}{Proof} We insert into  \thetag{\ref{fml:17}} the elements $a=X_1\wedge
\dots\wedge X_n$ and $\mu=du^1\wedge \dots\wedge du^m$. Since
\begin{equation*}
\mathcal L_{X_i}a=\sum_k c^k_{ik}a\quad\text{  and }\quad \mathcal
L_{X_i}\mu=\sum_a\frac{\p \rh^a_i}{\p u^a}\mu,
\end{equation*}
we get
\begin{displaymath}
\langle X_i,\ph \rangle a\otimes \mu =
\Bigl(\sum_{k}c^k_{ik}+\sum_{a}\frac{\p\rh^a_i}{\p u^a}\Bigr)a\otimes\mu.
\qed\end{displaymath}
\end{demo}

One could say that representing cohomology locally does not make much
sense, e.g. the modular class $\on{Mod}(TM)$ is trivial so locally
trivial. However, remember that for a general Lie algebroid the Poincar\'e
lemma does not hold: closed forms need not be locally exact. In
particular, for a Lie algebra (with structure constants),
\thetag{\ref{fml:18}} says that the modular class is just the trace of the
adjoint representation. In any case,  \thetag{\ref{fml:18}} gives us a
closed form, which is not obvious on first sight. If $E$ is a trivial
bundle, \thetag{\ref{fml:18}} gives us a globally defined modular class in
local coordinates.

\subsection{Remark}\label{nmb:4.5} 
As we have already mentioned, the modular class of a Lie algebroid is
the first characteristic class of R.~L.~Fernandes \cite{Fe}. There are
also higher classes, shown in \cite{Cr} to be characteristic classes of
the anchor map, interpreted as a representation "up to homotopy". It is
interesting if our idea can be adapted to describe these higher
characteristic classes as well.

\section{Acknowledgement }\label{nmb:5}

The authors are grateful to Marius Crainic and Yvette Kosmann-Schwarzbach
for providing helpful comments on the first version of this note.

\bibliographystyle{plain}

\end{document}